\documentclass[smallextended,nospthms]{svjour3}
\RequirePackage{fix-cm}
\usepackage{amssymb}
\usepackage{amsmath}
\usepackage{amsfonts}
\usepackage{color}
\usepackage[ruled,vlined]{algorithm2e}
\usepackage{hyperref}
\usepackage{graphicx}
\usepackage{graphics}
\usepackage[numbers]{natbib}
\usepackage{lineno}
\usepackage{comment}
\newtheorem{Proposition}{Proposition}[section]
\newtheorem{Theorem}{Theorem}[section]

\newtheorem{Corollary}{Corollary}[section]
\newtheorem{Example}{Example}[section]

\newcommand\mybox{\hbox to 0pt{}\hfill$\rlap{$\sqcap$}\sqcup$}

\def\qed{\hbox to 0pt{}\hfill$\rlap{$\sqcap$}\sqcup$}

\journalname{Preprint}

\begin{document}

\baselineskip = 1.5\baselineskip

\title{Approximate solutions to the multiple-choice  knapsack problem by multiobjectivization and Chebyshev scalarization}

\author{Ewa M. Bednarczuk \and Ignacy~Kaliszewski \and Janusz Miroforidis}

\institute{Ewa M. Bednarczuk 
	\at Warsaw University of Technology \\
	Faculty of Mathematics and Information Science \\
	ul. Koszykowa 75, 00-662 Warszawa, Poland \\
	\email{Ewa.Bednarczuk@pw.edu.pl} \\
	ORCID: 0000-0003-3683-6881 \\
	\and Ignacy Kaliszewski
	\at Systems Research Institute, Polish Academy of Sciences\\ 
	ul. Newelska 6, 01-447 Warszawa, Poland \\ \email{Ignacy.Kaliszewski@ibspan.waw.pl} \\
	Warsaw School of Information Technology \\
	ul. Newelska 6, 01-447 Warszawa, Poland \\
	ORCID: 0000-0001-5404-7400 \\
	\and Janusz Miroforidis (corresponding author)
	\at Systems Research Institute, Polish Academy of Sciences\\ 
	ul. Newelska 6, 01-447 Warszawa, Poland \\ \email{Janusz.Miroforidis@ibspan.waw.pl \\
	ORCID: 0000-0002-1319-6239} \\
}

	
\maketitle

\begin{abstract}
	
The method BISSA, proposed by Bednarczuk, Miroforidis, and Pyzel, provides approximate solutions to the multiple-choice knapsack problem. To fathom the optimality gap that is left by BISSA, we present a method that starts from the BISSA solution and it is able to provide a better approximation and in consequence a tighter optimality gap. Like BISSA, the new method is based on the multiobjectivization of the multiple-choice knapsack but instead of the linear scalarization used in BISSA, it makes use of the Chebyshev scalarization.

We validate the new method on the same set of problems used to validate BISSA.  

\end{abstract}

\section{Introduction}
\label{intro} 

The multi-dimensional knapsack problem and the
multiple-choice knapsack problem (MCKP) are the two most noted generalizations of the knapsack
problem (KP) and are applied to modeling many real-life problems, e.g., in
project (investments) portfolio selection \cite{Nauss1978},\cite{Sinha1979}, capital budgeting \cite{Pisinger2001}, advertising
\cite{Sinha1979}, component selection in IT systems \cite{Kwong2010},\cite{Pyzel2014}, computer networks management \cite{Lee1999}, adaptive multimedia systems \cite{Khan1998}, and other.

The multiple-choice knapsack problem is formulated as follows. Given
are $m$ sets (categories) $N_1, N_2,\dots, N_m$ of items, of cardinality $|N_j| = n_j , \ j = 1,\dots,m$.  Real-valued nonnegative ‘profit’ $p_{ij} \geq 0$ and ‘cost’
$c_{ij} \geq 0, \ j = 1,\dots,m, \ i = 1,\dots,n_j$, is assigned to each item of each set.
The problem consists in choosing exactly one item from each set $N_j$, so that the
total cost does not exceed a given $b > 0$ and the total profit is maximized.

Let $x_{ij}, \ j = 1,\dots,m, \ i = 1,\dots,n_j$\, be defined as
\[
\begin{array}{rl}
x_{ij} = & \left\{ \begin{tabular}{c} 
	     $1$ if item $i$ from set $N_j$ is chosen, \\ \\
	     $0$ otherwise.
	     \end{tabular}
         \right.
\end{array}
\]

Note that all $x_{ij}$ form a vector $x$ of length $n = \sum^m_{j=1} n_j , \ x \in \{ 0,1 \}^n$, where
$x := (x_{11}, \dots, x_{1n_1}, x_{2_1},\dots, x_{2n_2},\dots, x_{m_1}, x_{m_2},\dots, x_{mn_m})^T$.

The {\em multiple-choice knapsack problem} takes the form
\begin{equation}
\begin{array}{l}
	\label{i_eq1}
	\tag{MCKP}
	\begin{array}{l}\max\ \ \sum_{j=1}^{m}\sum_{i=1}^{n_{j}}p_{ij}x_{ij} \\ \\
		\mbox{subject to:} \\ \\
		\sum_{j=1}^{m}\sum_{i=1}^{n_{j}}c_{ij}x_{ij} \leq b, \\ \\
	x \in X:=\{(x_{ij})\ |\ \sum_{i=1}^{n_{j}}x_{ij}=1, \\ \\
	x_{ij}\in\{0,1\}, \ j=1,\dots,m, \ \ i=1,\dots,n_{j}\}.
\end{array}
\end{array}
\end{equation}

Elements $x \in X$ are {\em feasible solutions} to MCKP if they satisfy \[\sum_{j=1}^{m}\sum_{i=1}^{n_{j}}c_{ij}x_{ij} \leq b
, 
\]and are {\em infeasible} solutions, if otherwise.

Bednarczuk et al. (\cite{Bednarczuk_al_2018}) proposed the method BISSA to approximate the optimal solution to MCKP by the problem multiobjectivization and the linear scalarization. They achieved this by providing Pareto optimal solutions to the resulting bi-objective optimization problem and exploiting effectively the specific structure of set $X$. Namely, by the multiple-choice constraints $\sum_{i=1}^{n_{j}}x_{ij}=1$, the  objective function $\sum_{j=1}^{m}\sum_{i=1}^{n_{j}}p_{ij}x_{ij}$ and the constraint function $\sum_{j=1}^{m}\sum_{i=1}^{n_{i}}c_{ij}x_{ij}$ can be calculated in each 
$j$th
category independently. In other words, these two functions are additively separable. Thus, a linear scalarizing function over these two functions is also additively separable. Clearly, this observation generalizes to any number of functions structured analogously. The method BISSA provides an approximate solution to MCKP together with the optimality gap estimation. 

A conceivable improvement of the method presented in \cite{Bednarczuk_al_2018} would be replacing the linear scalarization, that can provide only a subset of Pareto optimal solutions, with the Chebyshev scalarization, that in the case of MCKP (it is a discrete optimization problem) can provide all Pareto optimal solutions (cf. \cite{Miettinen1999}, \cite{Ehrgott2005}, \cite{Kaliszewski1987}, \cite{Kaliszewski2006}, \cite{Kaliszewski_al_2016}). 
However, the Chebyshev scalarizing function, as it is easy to show by a counterexample (as that one given in Appendix, cf. also \cite{Brownlee_al_2020}), is not additively separable. Thus, with the Chebyshev scalarization the specific structure of set $X$ and the additive separability of functions $\sum_{j=1}^{m}\sum_{i=1}^{n_{j}}p_{ij}x_{ij}$, $\sum_{j=1}^{m}\sum_{i=1}^{n_{j}}c_{ij}x_{ij}$ cannot be directly exploited for solving MCKP. Nevertheless, as presented in the next section, the Chebyshev scalarization may provide better approximate solutions to MCKP than the linear scalarization.

In this work, we exploit the multiobjectivization of MCKP and the Chebyshev scalarization of the resulting bi-objective problem, for constructing a method that attempts to find better approximate solutions to MCKP, as compared to those that can be found by the method proposed in \cite{Bednarczuk_al_2018}. 
We validate the new method on the same set of problems used to validate BISSA. Despite that BISSA produces approximate solutions within tight optimality gaps, especially for correlated objective function and constraint coefficients, in 20\% of instances of test problems our method was able to produce better results than that produced by BISSA. 

The outline of the paper is as follows. In Section \ref{multiobjectivization}, we present a bi-objective formulation of MCKP  and its Chebyshev scalarization. In Section 3, we discuss Pareto optimality of additively separable functions. We exploit the Chebyshev scalarization in Section  \ref{section4} where we present a procedure to narrow the optimality gap to MCKP, as compared to the optimality gap provided by \cite{Bednarczuk_al_2018}, together with its effective implementation. Section \ref{experiments} reports results of numerical experiments. Section \ref{conclusions} concludes.

\section{Multiobjectivization of MCKP}
\label{multiobjectivization}

Given function $F: X \rightarrow \mathbb{R}^{k}$, $F = (f^1, \dots, f^k)$, $f^l : X \rightarrow \mathbb{R}, \ l = 1,\dots,k$, solution $\bar{x} \in X$  is Pareto optimal on $X$, if
$$
(F(\bar{x}) + \mathbb{R}^+_k) \cap F(X) = F(\bar{x})
$$
where $\mathbb{R}^+_k = \{ y  \in \mathbb{R}^k \, | \, y_l \geq 0, \ l = 1,\dots,k \}$. 

We multiobjectivize MCKP as follows:
we keep maximizing the original MCKP objective while minimizing the new objective originated from the constraint. For the sake of clarity of the presentation, in the sequel we replace  minimization of the second objective function by  maximization of its negative.
In this way, we obtain the following bi-objective (maximization) problem
\begin{equation}
	\label{i_eq2}
	\tag{MO\_MCKP}
\begin{array}{l}
	    \mbox{vmax} \ (f^1(x), f^2(x))  \\ \\
	   	\mbox{subject to:} \\ \\
		x \in X ,
\end{array}
\end{equation}
where $f^1(x) = \sum_{j=1}^{m}\sum_{i=1}^{n_{j}}p_{ij}x_{ij}$, $f^2(x) = -\sum_{j=1}^{m}\sum_{i=1}^{n_{j}}c_{ij}x_{ij}$, 'vmax' denotes the operator of deriving Pareto optimal solutions.

In this formulation, the structure of the multiple-choice set $X$ is fully exposed. The advantage of such a formulation is that now the feasible set of \eqref{i_eq2} is just $X$ and  not the set of $x \in X$ restrained by the additional linear inequality, as it is in the case of the original MCKP problem.

For a general reformulation of constrained singleobjective problems as multiobjective problems see \cite{Klamroth_Tind_2006}.

In Bednarczuk et al. (\cite{Bednarczuk_al_2018}), the search for approximate solutions to MCKP was conducted among solutions to MO\_MCKP that are Pareto optimal and derived by the linear scalarization (commonly referred to as {\it supported solutions}). Since, as observed above, a linear scalarizing function over two objective functions of MO\_MCKP is additively separable, such solutions are Pareto optimal for $F_j(x) = (f_j^1(x), f_j^2(x)) = (\sum_{i=1}^{n_{j}}p_{ij}x_{ij}, -\sum_{i=1}^{n_{j}}c_{ij}x_{ij})$ on $X_{j} :=\{ x_j = (x_{ij})\ |\ \sum_{i=1}^{n_{j}}x_{ij} = 1, \
x_{ij}\in\{0,1\}  \}$, $j = 1,\dots,m$, and Pareto optimal for $F(x) = (f^1(x), f^2(x))$ on $X$. Sets $X_j, \ j =1,\dots,m$, are referred to as components of $X$. The structure of MO\_MCKP and the adopted notation is represented in Table \ref{table1}.

\begin{table}
	\label{table1}
\[
	\begin{array}{|c|c|c|c|c|c|}
		\hline 
		F(x)  & F_1(x_1) & \dots & F_j(x_j) & \cdots & F_1(x_m) \\ \hline \hline 
		f^1(x)  & f^1_1(x_1) & \dots & f^1_j(x_j) & \cdots & f^1_1(x_m) \\ \hline
		f^2(x)  & f^2_1(x_1) & \dots & f^2_j(x_j) & \cdots & f^2_1(x_m) \\ \hline \hline
		 \multicolumn{1}{c|}{}& x_1 \in X_1 & \cdots & x_j \in X_j & \cdots & x_m \in X_m \\ \hline 
		X = & X_1 & \times \cdots \times & X_j & \times \cdots \times & X_m \\ \hline
	\end{array}
\]
	\caption{The structure of MO\_MCKP and the adopted notation.}
	\label{table1}
\end{table}

In Bednarczuk et al. (\cite{Bednarczuk_al_2018}) the following result was proved.
Let $S_P$ be the set of all Pareto optimal solutions to MCKP.

\begin{Theorem}
	\label{theorem2.1}
Let $x^* \in X$ be a Pareto optimal solution to (MO\_MCKP), such that
\[
	b - c^Tx^* = \min_{ \{x \in S_P, \ b - c^Tx \geq 0 \} } \ b - c^Tx. 
\]
Then $x^*$ solves (MCKP).
\end{Theorem}

The authors limited their considerations to the subset of $S_P$, namely to set $S_{\bar{P}}$ of the so-called supported Pareto optimal solutions, i.e. solutions that can be derived by linear scalarizing functions, a task that in the case of multiple-choice constraints reduces to sorting. In consequence, by algorithm BISSA  proposed they were able to establish the optimality gap with respect to $S_{\bar{P}}$, that amounts to 
\[
b - c^T\bar{x}^* = \min_{ \{x \in S_{\bar{P}}, \ b - c^Tx \geq 0 \} } \ b - c^Tx \,. 
\]
Actually, BISSA finds two solutions to MO\_MCKP: $x^*$ and $x^{'*}$, that solve the following problem 
\begin{equation}
	\label{gap}
c^Tx^* - c^Tx^{'*} = \min_{\, x \in \{S_{\bar{P}}, \ b - c^Tx \geq 0 \} \, , \, x' \in \{S_{\bar{P}}, \ b - c^Tx \leq 0 \} } \,
 c^Tx - c^Tx'\, .
\end{equation}

Our idea developed in this work is to narrow the optimality gap left by algorithm BISSA. The idea relies on exploiting  $x\in X$, such that $x_{j} \in X_{j}$, are Pareto optimal for $F_{j}$ on $X_{j}$ but not necessarily supported on $X_j$, $ j=1,\dots,k$.
 

For completeness, we formulate the Chebyshev scalarization for MO\_MCKP.

Let $y^{*1} > \max_{x \in X} \sum^m_{j=1}\sum^{n_j}_{i=1}p_{ij}x_{ij}, \ y^{*2} >  \max_{x \in X} -\sum^m_{j=1}\sum^{n_j}_{i=1}c_{ij}x_{ij}$.

The Chebyshev scalarization takes the following form 

\begin{equation}
	\label{i_eq3}
	\begin{array}{c}
		\min_{x \in X} \max \ \left\{ \begin{array}{c} \lambda_1(y^{*1} - f^1(x))  \\ \\ \lambda_2(y^{*2} - f^2(x)) 
		 \end{array} \right. \\ \\
		\mbox{subject to:} \\ \\
		x \in X.	
	\end{array} 
\end{equation}

As observed, this function is not additively separable.

\section{Pareto optimality of additively separable functions}
\label{Pareto_decomposable}

Below, we exploit the following general facts concerning additively separable functions. The following theorem states the underlying fact for the procedure presented in Section \ref{section4}.


\begin{Theorem}
	\label{theorem3_1}
	Given function $F: Y \rightarrow \mathbb{R}^{k}$, $F = (f^1, \dots, f^k)$,
	defined on $Y = Y_1 \times \cdots \times Y_m$, $Y_j \in \mathbb{R}^{n_j}, \ j =1,\dots,m$, each $f^l:Y\rightarrow \mathbb{R}$, $l=1,\dots,k$, is additively separable on $Y$, i.e., $f^l:=\sum_{j=1}^{m} f^{l}_{j}$, $l=1,\dots,k$,  where $f^{l}_{j}:Y_{j}\rightarrow\mathbb{R}$, $j=1,\dots,m$.
	
	Then, for every Pareto optimal solution $x = (x_1, \dots , x_m)$ of $F$ on $Y$, $x_j \in Y_j$ is a Pareto optimal solution of $F_j=(f_{j}^{1},\dots,f_{j}^{k})$ on $Y_j, \ j = 1,\dots,m$. 
\end{Theorem}

\noindent
\textbf{Proof}. Suppose that on the contrary, for given $\bar{x}=(\bar{x}_{1},\dots,\bar{x}_{m})$ there exists index $j$, $1 \le j \le m$,  such that $\bar{x}_{j} \in Y_{j}$ is not Pareto optimal of $F_{j}$  on $Y_j$. Then there exists $x_{j}\in Y_{j}$ such that
$f^l_{j}(x_{j}) \geq f^l_{j}(\bar{x}_{j}), \ l = 1,\dots,k$, and $f^l_j(x_j) > f^l_j(\bar{x}_j)$ for some $l$.  Hence, for $x' = (\bar{x}_1,\dots,x_j,\dots,\bar{x}_m)$, one has $f^l(x') \geq f^l(\bar{x}), \ l = 1,\dots,k$, and for some $l, \ f^l(x') > f^l(\bar{x}_j)$, which contradicts Pareto optimality of $\bar{x}$. \mybox

The converse statement does not hold, i.e. it is not true that $x$ composed of Pareto optimal solutions $x_{j}$ of $F_{j}$ on $Y_{j}, \ j=1,\dots,m$, \ $x=(x_{1},\dots,x_{m})$, is the Pareto optimal solution of $F$ on $Y$. (see the example in the Appendix, cf. also \cite{Brownlee_al_2020}). 

By Theorem \ref{theorem3_1},
searching for Pareto optimal solutions of $F$ on $Y$ can be limited to the subset of $Y$, such that for all $j =1,\dots,m$, $x_j$ is Pareto optimal of $F_{j}$ on $Y_j$.

\begin{Corollary}
	\label{corollary3_1}
	Given function $F: Y \rightarrow \mathbb{R}^{2}$, $F = (f^1, f^2)$,
	defined on $Y = Y_1 \times \dots \times Y_m$, each $f^i$, $i=1,2$, is additively separable with respect to $Y$.
	
	Let $x^A \in Y$ and $x^B \in Y$ be such that $f^1(x^A) < f^1(x^B)$.  
	
	If 	the set $\Phi:=\{ x \in Y \ | \ f^1(x) > f^1(x^A), \ f^2(x) > f^2(x^B) \}$ contains no  element $(x_{1},\dots, x_{m})$ such that $x_{j}$ is Pareto optimal of $F_{j}$ on $Y_j$ for each $j=1,\dots,m$, then $\Phi$ contains no Pareto optimal solution of $F$ on $Y$.
\end{Corollary}

\noindent
\textbf{Proof}. The assertion of this corollary follows immediately from Theorem \ref{theorem3_1}. \mybox

\section{A procedure to narrow the optimality gap to MCKP}
\label{section4}

We now deal with $F = (f^1, f^2)$ where $f^1(x) = \sum^m_{j=1}\sum^{n_j}_{i=1}p_{ij}x_{ij}$, $f^2(x) = -\sum^m_{j=1}\sum^{n_j}_{i=1}c_{ij}x_{ij}$.

The procedure to narrow the optimality gap to MCKP, as compared to that provided by \cite{Bednarczuk_al_2018},  is based on Theorem \ref{theorem3_1}.

\subsection{The procedure}
\label{procedure}

The procedure starts where algorithm BISSA stops, namely it starts with two Pareto optimal solutions to MCKP that solve problem (\ref{gap}): $x^*$ that is feasible to MCKP (below denoted $x^A$) and $x^B$ that is  infeasible to MCKP (below denoted $x^B$). BISSA guarantees that above the line passing through $f(x^A)$ and $f(x^B)$ there is no $f(x)$ such that $x$ is Pareto optimal solution. Thus, the set that can contain better feasible solutions to MCKP that $x^A$ is defined by $ f_1(x^A) \leq f_1(x^B)$, $ f_2(x^A) \geq f_2(x^B)$, $f(x) \leq_{\mathbb{R}^2_+} \lambda f(x^A) + (1 - \lambda) f(x^B), \ 0 < \lambda < 1$.

Let $\varepsilon>0$. To simplify the presentation, we assume for a while that the Chebyshev scalarization always yields Pareto optimal solutions, whereas in general it provides weakly Pareto optimal solutions. We will show how to tackle this problem in the next section. With this assumption the algorithm is as follows.

Let $J^{1}$ be the set of all $1\le j\le m$ such that
$f^1_j(x^A_j) < f^1_j(x^B_j)$, i.e.,
\begin{equation} 
	\label{jot1}
J^{1}:=\{1\le j\le m\ | \ f^1_j(x^A_j) < f^1_j(x^B_j), \ f^2_j(x^A_j) > f^2_j(x^B_j) \}.
\end{equation}

Observe that at the starting iteration, $J^{1}\neq\emptyset$. Indeed, if $J^1$ were empty, then for all $j = 1,\dots,m$, 
$$
f^1_j(x^A_j) < f^1_j(x^B_j),$$
and
$$ 
f^2_j(x^A_j) \leq f^2_j(x^B_j),
$$
hence $x^A$ would be dominated by $x^B$, a contradiction to Pareto optimality of $x^A$.


Observe also that since on the first iteration $x^{A}$ is a Pareto optimal solution for $(f^{1},f^{2})$ on $X$, then, by Theorem \ref{theorem3_1}, each $x^{A}_{j}$ is a Pareto optimal solution of $F_{j}$ on $X_j, \ j=1,\dots,m$. Hence, for each $j\in J^{1}$, we have
\begin{equation}
	\label{eq_ineq}
f^{2}_{j}(x^{A}_{j})>f^{2}_{j}(x^{B}_{j}).
\end{equation}
This  is illustrated in Figure 1.

\begin{figure}
	\centering
	\includegraphics[scale=0.5]{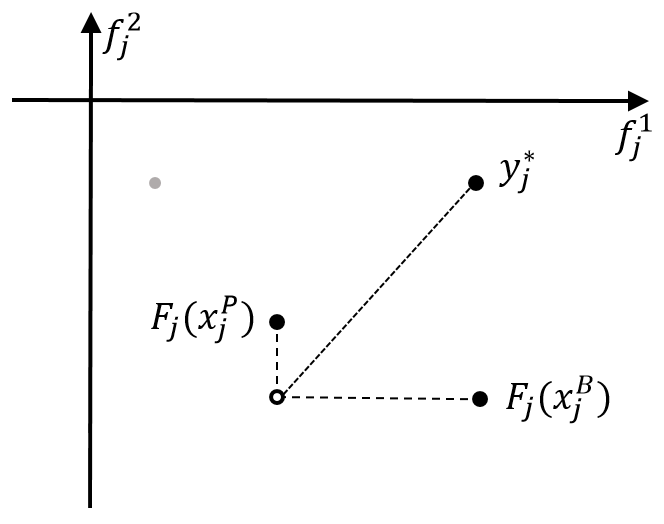}
	\caption{''The layout'' at the first iteration of the algorithm for some $j$. Element $\circ$ defines weights of the Chebyshev scalarization.}
	\label{Fig1}
\end{figure}


\vspace{0.3cm}
\noindent{\it $s$-th iteration of the procedure, $s = 1,\dots, \, .$ Each iteration consists of two steps.}

\noindent \textbf{Step 1} 

Determine $J^{1}$. 

For each $j\in J^{1}$, 
determine a Pareto optimal solution $x^P_j$ by solving the Chebyshev problem
\[
\min_{x_j \in X_j} \max_{i = 1,2} \ \lambda^j_i(y^{*i}_j - f^i_j(x_j))
\]
where: 
\[
y^{*1}_j = \max_{x_j \in X_j} f^1_j(x_j) + \varepsilon,
\]
\[
y^{*2}_j = \max_{x_j \in X_j} f^2_j(x_j) + \varepsilon,
\]
\[
\lambda^j_1 = (y^{*1}_j - f^1_j(x^A_j))^{-1},
\]
\[
\lambda^j_2 = (y^{*2}_j - f^2_j(x^B_j))^{-1}.
\]  

Let $J \subseteq J^1$ be such that for $j \in J$, $f^1_j(x^A_j) < f^1_j(x^P_j)$. 

Observe that for each $j \in J$, we have 
$$ f^1(x^A) < f^1(x^{A_j}), \ \mbox{where} \ x^{A_j} = (x^A_1, \dots, x^P_j, \dots, x^A_m).
$$ 

\noindent \textbf{Step 2}

Let $J^b, \ J^b \subseteq J$, denote the set of indices $j$ such that $x^{A_j}$ are feasible, i.e. $f^2(x^{A_j}) \geq -b$.


\noindent
Choose any $j^* \in J^b$. Set $x^A := x^{A_{j^*}}$ and proceed to the next iteration.

\vspace{0.3cm}

\begin{Proposition}
	If $J$ is empty and $x^A$ is a Pareto optimal solution of \eqref{i_eq2}, then $x^A$ is the optimal solution to \ref{i_eq1}.
\end{Proposition}

\noindent
\textbf{Proof}. If $J$ is empty, then by Corollary \ref{corollary3_1} the set
\[
\{x \in X | \ f^1_j(x^A_j) \leq f^1_j(x^B_j), \ f^2_j(x^A_j) \geq f^2_j(x^B_j) \}
\]
contains only two Pareto optimal solutions to MO\_MCKP, $x^A$ that is feasible to MCKP, and $x^B$ that is infeasible to MCKP, thus $x^A$ is optimal to MCKP.
\mybox

\begin{Proposition}
	If $J^b$ is empty and $x^A$ is a Pareto optimal solution of \eqref{i_eq2}, then $x^A$ is the optimal solution to \eqref{i_eq1}.
\end{Proposition}

\noindent
	\textbf{Proof}. If $J^b$ is empty, then by Corollary \ref{corollary3_1} the set
	\[
	\{x \in X | \ f^1_j(x^A_j) \leq f^1_j(x^B_j), \ f^2_j(x^A_j) \geq f^2_j(x^B_j) \}
	\]
	contains only one Pareto optimal solution to MO\_MCKP, namely $x^A$, that is feasible to MCKP, thus $x^A$ is optimal to MCKP.
\mybox

Observe that starting from the second iteration there is no guarantee that $x^{A}$ is a Pareto optimal solution for $(f^{1},f^{2})$ on $X$. Hence, from the second iteration on, it can happen that
\begin{equation}
	\label{eq_4.3}
	f^{2}_{j}(x^{A}_{j})<f^{2}_{j}(x^{B}_{j})\ \ \text{for some  } j\in J^{1}.
\end{equation} 

\vspace{0.7cm}
Different rules for selecting $j^* \in J^{b}$ and thus $x^{P}_{j}$ 
can be proposed, for example $\tilde{x}^{P}_{j^*}$ $ = \arg\max_{j \in J^b} f^1(x^{P}_{j})$ (the MAX-PROFIT rule).

\begin{figure}
	\centering
	\includegraphics[scale=0.5]{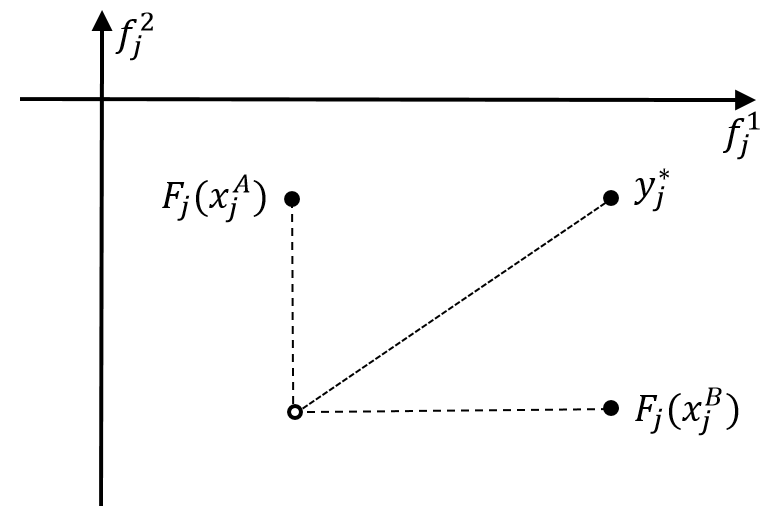}
	\caption{$s$-th iteration of the algorithm, $s > 1$, for some $j^b$ in the case the Pareto optimal solution $x^A$, derived by the Chebyshev scalarization in iteration $s-1$, is feasible; $x^P_j$, replaces $x^A_j$ from the previous iteration (now in the light gray color); $\circ$ is the element defining new weights in the Chebyshev scalarization.} 
	\label{Fig2}
\end{figure}

Figure \ref{Fig2} illustrates operations of one iteration of the procedure for some $j^*$.

The procedure stops when either set $J$ or set $J^b$ becomes empty.

\begin{algorithm}
	\SetAlgorithmName{Algorithm}{} 
	\scriptsize
	\caption{Procedure KISSA to solve MCKP.} 
	\label{proc_1} 
	\vspace{0.2cm}
	\nl \textbf{INPUT}: \\ $x^A = (x^A_1,\dots,x^A_m)$ -- a feasible solution to MCKP such that $x^A_j$ is Pareto optimal to $F_j$ on $X_j, \ j=1,\dots,m$, provided by BISSA;\\
	\vspace{0.2cm}
	 $x^B = (x^B_1,\dots,x^B_m)$ -- an infeasible solution to MCKP such that $x^B_j$ is Pareto optimal to $F_j$ on $X_j, \ j=1,\dots,m$, provided by BISSA ($f^1(x^A) < f^1(x^B)$). \\
	\vspace{0.2cm}
	\nl \textbf{OUTPUT}: a feasible solution to MCKP with a greater than or equal value of the objective function for the input feasible solution $x^A$. \\
	\vspace{0.2cm}
	\nl \Begin{ \vspace{0.2cm}
		\nl Determine $\rho$.\\ \vspace{0.2cm}
		\nl \While{\textbf{true}}
		{
			\vspace{0.2cm}
			\nl Determine $J^1$. \\
			\vspace{0.2cm}
			\nl $J:=\emptyset$. \\
			\vspace{0.2cm}
			\nl \For{\textbf{each} $j \in J^1$}
			{	
					\vspace{0.2cm}
				\nl Determine $y^{*}_j = (y^{*1}_j, y^{*2}_j )$ by setting: \\
					\vspace{0.2cm}
				\nl \hspace{0.5cm} $y^{*1}_j := \max_{x_j \in X_j} f^1_j(x_j) + \varepsilon$, \\
					\vspace{0.2cm}
				\nl \hspace{0.5cm} $y^{*2}_j := \max_{x_j \in X_j} f^2_j(x_j) + \varepsilon$.	  \\
	\vspace{0.2cm}
				\nl $\lambda^j_1 := (y^{*1}_j - f^1_j(x^{A_j}))^{-1}$. \\
				\vspace{0.2cm}
				\nl $\lambda^j_2 := (y^{*2}_j - f^2_j(x^{B_j}))^{-1}$. \\
				\vspace{0.2cm}
				\nl Solve \ $\min_{x_j \in X_j} \max_{i = 1, 2} \ \lambda^i_j ((y^{*i}_j - f^i_j(x_j)) + \rho \sum^{m}_{i=1} (y^{*i}_j - f^i_j(x_j))$ \\
				\vspace{0.1cm} to derive $x^P_j$. \\
					\vspace{0.2cm}
				\nl \If{$f^1_j(x^A_j) < f^1_j(x^P_j)$}{\vspace{0.1cm} add $j$ to $J$.}
					\vspace{0.2cm}
			}
			\nl \If{$J = \emptyset$}{\vspace{0.1cm} \textbf{break}}
			\vspace{0.2cm}
			\nl Determine $J^b \subseteq J$. \\
			\vspace{0.2cm}
			\nl \If{$J^b = \emptyset$}{\vspace{0.1cm} \textbf{break}}
			\vspace{0.1cm}
				\vspace{0.2cm}
			\nl Select $j^* \in J^b$. \\
				\vspace{0.2cm}
			\nl $x^A := (x^A_1, \dots, x^P_{j^*}, \dots, x^A_m)$.
		}
	\vspace{0.2cm}
	\nl \textbf{RETURN $x^A$}. \\
	}
\end{algorithm}
The procedure to narrow the optimality gap to MCKP described above we shall call KISSA.

\subsection{The implementation of KISSA}
\label{algorithm}

Now, we have to correct for a temporary assumption  we have made that solutions to the Chebyshev problem, as above, are Pareto optimal, whereas they are in fact only weakly Pareto optimal. Pareto optimal solutions to $X_j, \ j = 1,\dots,m$, can be derived by the following result (\cite{Kaliszewski1987}, Theorem 4.6, p. 54) adapted to the notation of this work.

Denote
\[
\gamma = \{t \ | \ x^{t}_j \in X_j \,, \ x^t_j \mbox{ is Pareto optimal on $X_j$} \}\, .
\]

\begin{Theorem}
	\label{theorem4_1}
	Let
	\begin{equation}
		\label{eq_4.2}
	\rho < \min_{t\in \gamma}
	\left\{ \begin{array}{c}
		\min_{u\in\gamma \setminus\{t\}}
		\left\{ \begin{array}{c}
			\frac{\min_{l;(f^l_j(x^t_j) - f^l_j(x^u_j) > 0}(f^l_j(x^t_j) - f^l_j(x^u_j))}{\sum^{k}_{l=1}(f^l_j(x^u_j) - f^l_j(x^t_j))}  
		        \end{array} \right. 	         
	       \end{array} \right.
   \end{equation}
\[
  \hspace{5cm}          \left| \left. \sum^{k}_{l=1}(f^l_j(x^u_j) - f^l_j(x^t_j))>0 \right\} \right\} \, .
\]	
	
	Solution 		
	$\bar{x}_j \in X_j$ is Pareto optimal on $X_j$ if and only if there exists a vector
	$\lambda\, , \ \lambda>0 \, ,$ such that
	$\bar{x}_j$ solves 
	\[
	\min_{x_j \in X_j} \max_{1 \leq l \leq k} \ \lambda_l ((y^{*j}_l - f^l_j(x_j)) + \rho \sum^{k}_{s=1} (y^{*l}_s - f^s_j(x_j)),  
	\]
	\noindent
	where $y^{*j}_l > \max_{x_j \in X_j} f^l_j(x_j), \ l =1,\dots,k$.
\end{Theorem}

Let us denote the set of all Pareto optimal solutions by $N$.

To make the formula for $\rho$, given in Theorem \ref{theorem4_1}, operational, we can use a conservative upper bound on $\rho$, namely $\delta$, where

\[
\delta = \min_{x_j \in X_j}
\left\{
       \begin{array}{c}
	\min_{x'_j \in \ X_j \setminus\{x_j \}}
	\left\{ \begin{array}{c}
		\frac{\min_{l;(f^l_j(x_j) - f^l_j(x'_j) > 0}((f^l_j(x_j) - f^l_j(x'_j))}{\sum^{k}_{l=1}(f^l_j(x'_j) - f^l_j(x_j))} 
		    \end{array} \right.
	  \end{array} \right.
\]
\[
\hspace{5cm}          \left| \left. \sum^{k}_{l=1}(f^l_j(x_j) - f^l_j(x'_j))>0 \right\} \right\} \, .
\]	
In the case of MO\_MCKP, function values are just values of single coefficients. Therefore, the above formula takes the form 
\[
\delta = \min_{i, i=1,\dots,n_j}
\left\{
\begin{array}{c}
	\min_{i', i'=1,\dots,n_j, i' \neq i}
	\left\{ \begin{array}{c}
		\frac{\min_{l;(d^l_{ij} - d^l_{i'j}) > 0}((d^l_{ij} - d^l_{i'j})}{\sum^{k}_{l=1}(d^l_{ij} - d^l_{i'j})} 
	\end{array} \right.
\end{array} \right.
\]
\[
\hspace{5cm}          \left| \left. \sum^{k}_{l=1}(d^l_{ij} - d^l_{i'j}) > 0 \right\} \right\} \, .
\]	
where $d^l_{ij} = p^l_{ij}$ if $l = 1$, and $d^l_{ij} = -c^l_{ij}$ if $l = 2$.

\section{Numerical experiments}
\label{experiments}

We conducted numerical experiments on the same set of randomly generated test problems as in Bednarczuk et al. (\cite{Bednarczuk_al_2018}). We considered three sets of 10 problems each with uncorrelated coefficients of the objective function and the constraint, with, respectively, 10 categories and 1000 variables, 100 categories and 100 variables, 1000 categories and 10 variables. We also considered one set of 10 problems with weakly correlated coefficients, notoriously hard to solve for BISSA, with 20 categories and 20 variables. The results are given in Table \ref{Table1} and Table \ref{Table2}.
The second column of the tables is the relative gap between the objective function values of the solutions obtained by the exact algorithm used in \cite{Bednarczuk_al_2018} and BISSA.
The third column of the tables is the relative gap between the objective function values of the solutions obtained by the exact algorithm used in \cite{Bednarczuk_al_2018} and KISSA.
The last column of the table shows the number of improvements in the objective function value during KISSA's operation. If empty, no improvement is observed. In the numerical experiments, we used the following parameter values: $\rho = 1\text{E}-7, \ \varepsilon = 1\text{E}-4$.
To select index $j^{*}$ (line 19 of Algorithm 1), we applied the MAX-PROFIT rule shown in subsection \ref{procedure}.

The results show that compared with BISSA, which provides very tight optimality gaps and often optimal solutions, the room for improvements offered by KISSA is limited. Despite that, in 20\% of instances of test problems, our method was able to produce better results than that produced by BISSA. Significant tightening of the optimality gaps can be observed for problems with weakly correlated coefficients, where the optimality gaps provided by BISSA are much looser than in the case of uncorrelated coefficients.

As the computation time of KISSA for each of the problems tested was fractions of a second, it can be used as a no-cost plug-in for BISSA.
\begin{table}
	\centering
	\includegraphics[scale=0.65]{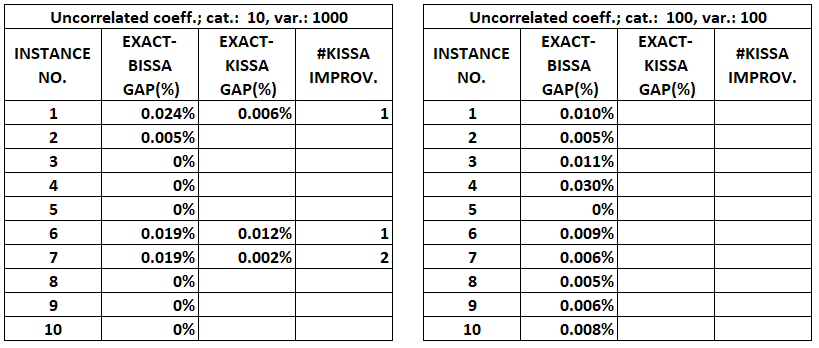}
	\caption{KISSA's results for uncorrelated instances.} 
	\label{Table1}
\end{table}

\begin{table}
	\centering
	\includegraphics[scale=0.65]{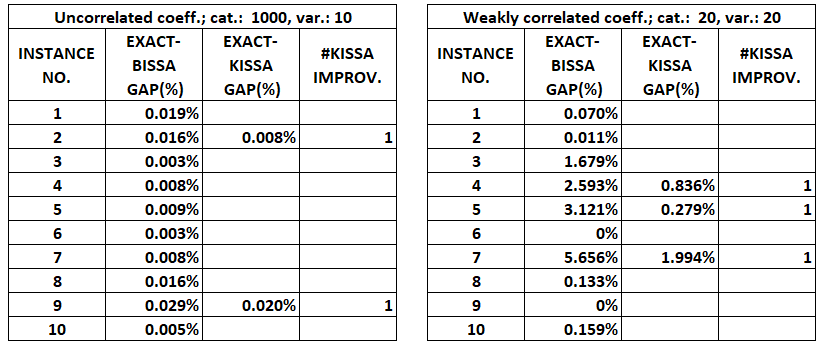}
	\caption{KISSA's results for uncorrelated and weakly correlated instances.} 
	\label{Table2}
\end{table}
\section{Conclusions}
\label{conclusions}

In view of the very good performance of BISSA, the extra effort required by KISSA is worth the trouble in problems where data are precise and the objective value reflects a practical problem of high stakes consequences. In such cases, a quest for approximate solutions tightly close to the optimal one is justified, and KISSA is a viable option. The decision of whether it is the case depends on the nature of practical problems modeled as MCKP.


\section*{Declarations}
\textbf{Conflict of interest} The authors declare that they have no conflict of interest.\\
\textbf{Funding} The authors did not receive support from any organization for the submitted work.

\section*{Appendix}

\begin{Example}
	\label{example}
	Assume $k =2, \ m = 2, \ n_1 = 2, \ n_2 =2$, 
	\[
	\left\{ 
	\begin{array}{rl}
	& p^l_{ij} \\
	- & c^l_{ij} 
	\end{array} 
    \right\}
	= \left\{
	\begin{array}{cccc}
		2    & 3  & 4  & 2 \\
		-1.9 & -3 & -2 & -1
	\end{array}
	\right\} .
	\]

The values of the objective functions for all four feasible  solutions $X = \{1010, 1001, 0110, 0101\}$	are given in Table \ref{table3}.

\begin{table}
\begin{tabular}{|c|c|c|c|c|}
	\hline
	& \multicolumn{4}{c|}{Solutions}    \\
	& Pareto \ opt. & 	Pareto \ opt. & 	Pareto \ opt. & dominated \\ \hline
	&  1\,0\,1\,0  & 1\,0\,0\,1 & 0\,1\,1\,0  &	0\,1\,0\,1    \\ \hline 
	$i=1$ & 6    & 4    & 7     &  -5   \\
	$i=2$ & -3.9  &-2.9	& -5 	&  -4  \\
	\hline
\end{tabular}
\caption{The values of the objective functions for all four feasible solutions (Example \ref{example}).} 
\label{table3}
\end{table}

Now, we derive Pareto optimal solutions in each component separately. In the first component they are  $(1\,0), \ (0\,1)$, and in the second component they are $(10), \ (01)$, as shown in Table \ref{table4}. 		
\begin{table}
\begin{tabular}{|c|c|c||c|c|}
		\hline
	& \multicolumn{4}{c|}{Solutions \ componentwise}   \\
	& \multicolumn{2}{c||}{$j=1$} & \multicolumn{2}{c|}{$j=2$} \\ \hline
	& Pareto \ opt. & Pareto \ opt. & Pareto \ opt. & Pareto \ opt. \\ \hline
	& 1\,0 & 0\,1 & 1\,0  & 0\,1 \\ \hline 
	$i=1$	& 2  & 3   & 7    & 5  \\
	$i=2$	& -1.9  & -3 & 4  & 2 \\
		\hline
\end{tabular}
\caption{Pareto optimal solutions in each component separately (Example \ref{example}).} 
\label{table4}
\end{table}

The concatenation of the Pareto optimal solution in the component $X_1$, $(0\,1)$, with the Pareto optimal solution in component $X_2$, $(0\,1)$, resulting in a feasible solution $(0\,1\,0\,1)$, is not Pareto optimal in $X = X_1 \times X_2$.
\end{Example}


\end{document}